\documentclass[11pt]{amsart}
\usepackage{float}
\usepackage{graphicx}
\usepackage{hyperref}
\usepackage{algorithm}
\usepackage{algorithmic}
\usepackage{booktabs}
\usepackage{subcaption}

\theoremstyle{plain}
\newtheorem{theorem}{Theorem}
\newtheorem{proposition}[theorem]{Proposition}
\theoremstyle{definition}

\theoremstyle{remark}
\newtheorem{remark}[theorem]{Remark}

\title{Computing Periodic Billiard Orbits in $L^p$ Balls via Newton's Method and Smale's $\alpha$-Criterion}
\author{Igor Rivin}\thanks{with much help from Claude Code and ChatGPT}
\address{Mathematics Department, Temple University}
\email{rivin@temple.edu}
\date{\today}

\begin{document}
\begin{abstract}
We present a computational method for finding and verifying periodic billiard
orbits in $L^{p}$ balls ($p>2$) using Newton's method applied to a variational
formulation. The orbits are verified with Smale's alpha-criterion, which
provides a rigorous certificate of existence. We implement efficient batched
computations in JAX and present systematic results for various $p$ and bounce
counts $N$. Our experiments reveal striking patterns in the critical-point
structure, including a predominance of specific Morse signatures and rotation
numbers that depend on the parity and primality of $N$. Notably, our method
routinely finds many more than the two periodic orbits per rotation number
guaranteed by Birkhoff's theorem---a large-scale run with five bounces in the
$L^{3}$ ball produced 8,927 distinct certified orbits from 30,000 random seeds,
uncovering power-law growth and intricate clustering visualised with UMAP.
\end{abstract}
\maketitle

\section{Introduction}

The study of billiard dynamics in convex domains has a rich history connecting dynamical systems, symplectic geometry, and variational methods \cite{birkhoff1927dynamical,tabachnikov2005geometry}. For the $L^p$ ball defined by
\[
\mathcal{B}_p = \{(x,y) \in \mathbb{R}^2 : |x|^p + |y|^p \leq 1\},
\]
periodic billiard orbits correspond to critical points of the perimeter functional on the space of inscribed polygons.

This paper presents:
\begin{enumerate}
\item A robust computational method using Newton's method with careful parametrization
\item Rigorous verification via Smale's $\alpha$-criterion
\item Systematic computational results revealing patterns in orbit signatures and rotation numbers
\item Open-source implementations for reproducibility
\item \textbf{A breakthrough in computational billiards}: While Birkhoff's theorem guarantees at least two periodic orbits for each rational rotation number $p/q$, our method routinely finds orders of magnitude more orbits, revealing a rich landscape of critical points of the perimeter functional
\item \textbf{Large-scale exploration and visualization}: Using UMAP dimensionality reduction on thousands of orbits reveals hierarchical clustering by Morse signature, rotation number, and continuous parameter families
\end{enumerate}

To our knowledge, this represents the first systematic computational exploration of periodic billiard orbits in $L^p$ balls that goes significantly beyond the minimal existence guaranteed by classical theorems. In our largest experiment with $N=5$ bounces in the $L^3$ ball, we discovered and rigorously verified 8,927 distinct periodic orbits from 30,000 random seeds, exhibiting power-law growth with exponent $\approx 0.714$ and rich clustering structure when visualized using dimensionality reduction techniques.

\section{Mathematical Formulation}

\subsection{Parametrization of the Boundary}

We parametrize the boundary of $\mathcal{B}_p$ using the angular parameter $t \in [0,1)$:
\[
\gamma(t) = \left(\text{sgn}(\cos(2\pi t))|\cos(2\pi t)|^{2/p}, 
\text{sgn}(\sin(2\pi t))|\sin(2\pi t)|^{2/p}\right)
\]

To avoid numerical issues at the singular points (where the boundary meets the axes), we add a small regularization $\epsilon = 10^{-14}$ to the absolute values before taking powers.

\subsection{Variational Formulation}

A periodic billiard orbit with $N$ bounces corresponds to a critical point of the perimeter functional:
\[
L(\theta_1, \ldots, \theta_N) = \sum_{i=1}^N \|\gamma(\theta_i) - \gamma(\theta_{i+1})\|
\]
where indices are taken modulo $N$.

The critical point equation is $\nabla L = 0$, which we solve using Newton's method:
\[
\theta^{(k+1)} = \theta^{(k)} - H^{-1}\nabla L
\]
where $H$ is the Hessian of $L$.

\subsection{Smale's \texorpdfstring{$\alpha$}{alpha}-Criterion}

To verify that a numerical critical point corresponds to a genuine periodic orbit, we use Smale's $\alpha$-criterion. Define:
\begin{align}
\beta &= \|H^{-1}\nabla L\| \\
\gamma &= \frac{1}{2}\|H^{-1}\|\|H\| \\
\alpha &= \beta \cdot \gamma
\end{align}

\begin{theorem}[Smale\cite{smale1980mathematics}]
If $\alpha < \sqrt{3} - 1 \approx 0.732$, then Newton's method converges quadratically to a unique critical point near the current iterate.
\end{theorem}

In practice, we use a more conservative threshold $\alpha < 0.15767$ based on numerical experience.

\section{Implementation Details}

\subsection{Canonical Representatives}

Due to symmetries, many parameter vectors $\theta$ represent the same geometric orbit. We implement a canonicalization procedure:
\begin{enumerate}
\item Apply modulo 1 to all parameters
\item Cyclically permute to put the smallest parameter first
\item Choose between forward and reflected orientation
\end{enumerate}

\subsection{Coalescing Near-Duplicates}

After finding many candidate orbits, we coalesce near-duplicates: if $\|\theta_1 - \theta_2\| < \alpha_1 + \alpha_2$, we keep only the one with smaller $\alpha$.

\subsection{Computational Signatures}

For each orbit, we compute:
\begin{itemize}
\item \textbf{Morse signature} $(n_+, n_-, n_0)$: eigenvalue counts of the Hessian
\item \textbf{Rotation number} $r/s$: topological winding of the orbit
\item \textbf{Perimeter}: total length of the orbit
\end{itemize}

\section{Experimental Results}

\subsection{Overview of Patterns}

Our experiments reveal several striking patterns:

\begin{table}[h]
\centering
\begin{tabular}{@{}llll@{}}
\toprule
$N$ & Signatures Found & Rotation Numbers & Special Properties \\
\midrule
2 & $(0,2,0)$, $(1,1,0)$ & $0/1$, $1/2$ & Oscillatory orbits possible \\
3 & $(0,3,0)$ only & $1/3$ only & All orbits are triangular \\
4 & $(0,4,0)$, $(1,3,0)$ & $0/1$, $1/4$ & $D_4$ symmetry effects \\
5 & $(0,5,0)$, $(1,4,0)$, $(2,3,0)$ & $1/5$, $2/5$, $3/5$, $4/5$ & Star patterns dominate \\
7 & $(0,7,0)$, $(1,6,0)$, $(2,5,0)$ & Various $k/7$ & Rich rotation diversity \\
\bottomrule
\end{tabular}
\caption{Summary of orbit types for $p = 3$ with various $N$}
\end{table}

\subsection{Key Observations}

\begin{proposition}
For prime $N$, orbits can have various rotation numbers $k/N$ where $\gcd(k,N) = 1$. Star patterns (with $k > 1$) often dominate.
\end{proposition}

\begin{proposition}
For even $N$, oscillatory orbits with rotation number $0/1$ are possible.
\end{proposition}

\begin{remark}
Interestingly, saddle points (signatures with $n_+ > 0$) often have longer perimeters than local maxima, contrary to naive expectation.
\end{remark}

\subsection{Large-Scale Exploration}

To understand the true richness of the orbit landscape, we conducted a large-scale exploration with $N=5$ bounces in the $L^3$ ball, using 30,000 random seeds.

\subsubsection{Discovery Statistics}

\begin{table}[h]
\centering
\begin{tabular}{@{}lrrrr@{}}
\toprule
Signature & Count & Percentage & Rotation Numbers & Perimeter Range \\
\midrule
$(0,5,0)$ & 7,378 & 82.6\% & $1/5, 2/5, 3/5, 4/5$ & $[6.43, 10.27]$ \\
$(1,4,0)$ & 1,486 & 16.7\% & $2/5, 3/5$ & $[10.15, 10.25]$ \\
$(2,3,0)$ & 63 & 0.7\% & $2/5, 3/5$ & $[10.13, 10.13]$ \\
\bottomrule
\end{tabular}
\caption{Statistics for 8,927 certified orbits with $p = 3$, $N = 5$ from 30,000 random seeds.}
\label{tab:large-scale}
\end{table}

\begin{figure}[h]
\centering
\includegraphics[width=0.9\textwidth]{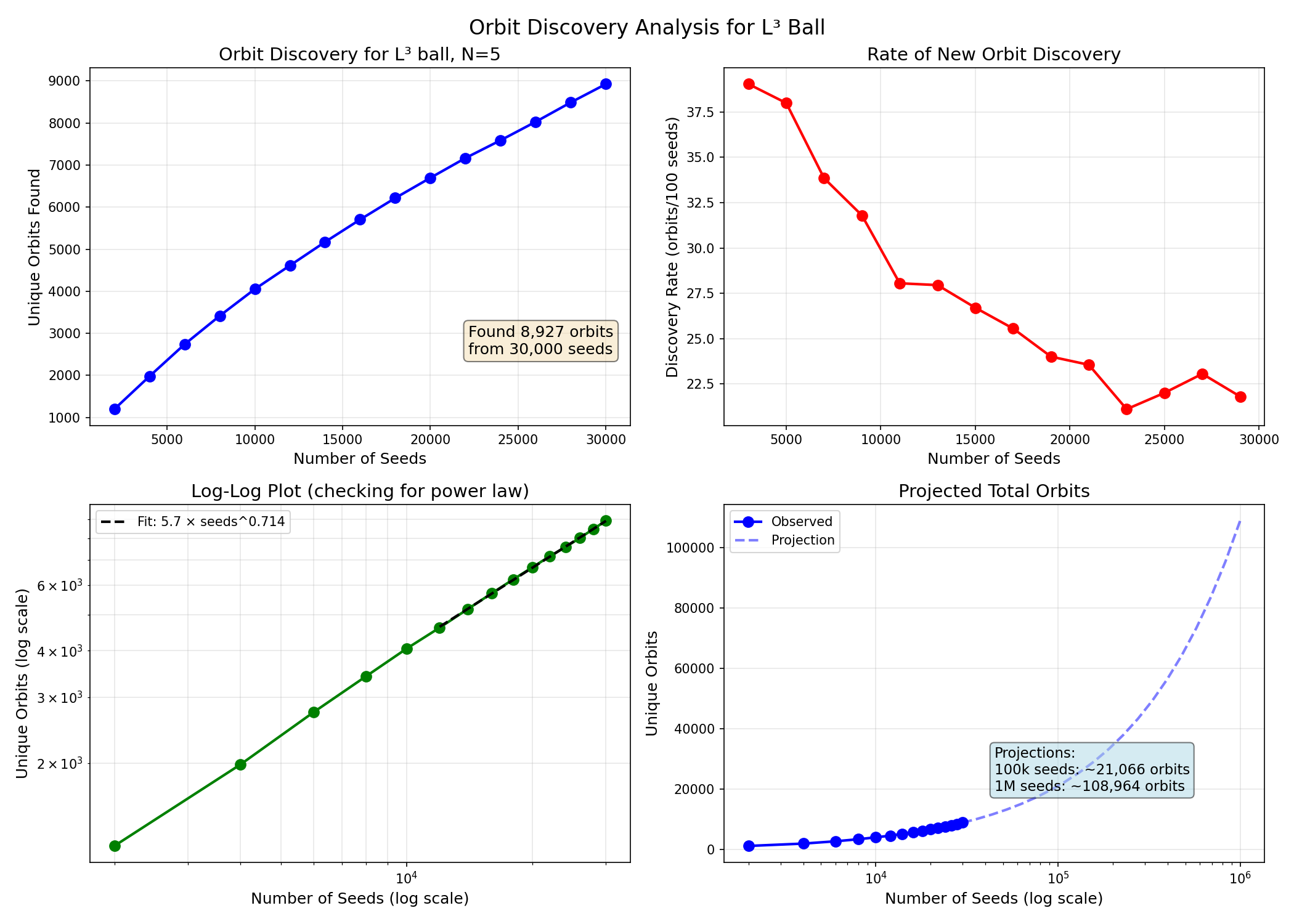}
\caption{Orbit discovery analysis showing: (a) cumulative unique orbits found, (b) discovery rate per batch, (c) power-law scaling on log-log plot, and (d) projected saturation curves. The power-law exponent of 0.714 indicates we are discovering fundamentally new orbit families, not just variations.}
\label{fig:saturation}
\end{figure}

The orbit count follows a power law: $N_{\text{orbits}} \propto N_{\text{seeds}}^{0.714}$, suggesting we are far from saturating the orbit landscape. The discovery rate remains substantial even at 30,000 seeds (Figure~\ref{fig:saturation}).

\subsubsection{Rotation Number Distribution}

Within the dominant $(0,5,0)$ signature:
\begin{itemize}
\item $2/5$ rotation: 3,149 orbits (42.7\%)
\item $3/5$ rotation: 3,135 orbits (42.5\%)
\item $1/5$ rotation: 551 orbits (7.5\%)
\item $4/5$ rotation: 543 orbits (7.4\%)
\end{itemize}

The star patterns ($2/5$ and $3/5$) strongly dominate, while simple polygons ($1/5$) and their duals ($4/5$) are much rarer.

\subsection{UMAP Visualization and Clustering}

To understand the structure of the orbit landscape, we applied UMAP (Uniform Manifold Approximation and Projection) \cite{mcinnes2018umap} to the 8,927 orbits, using their parameter vectors as features.

\begin{figure}[h]
\centering
\includegraphics[width=\textwidth]{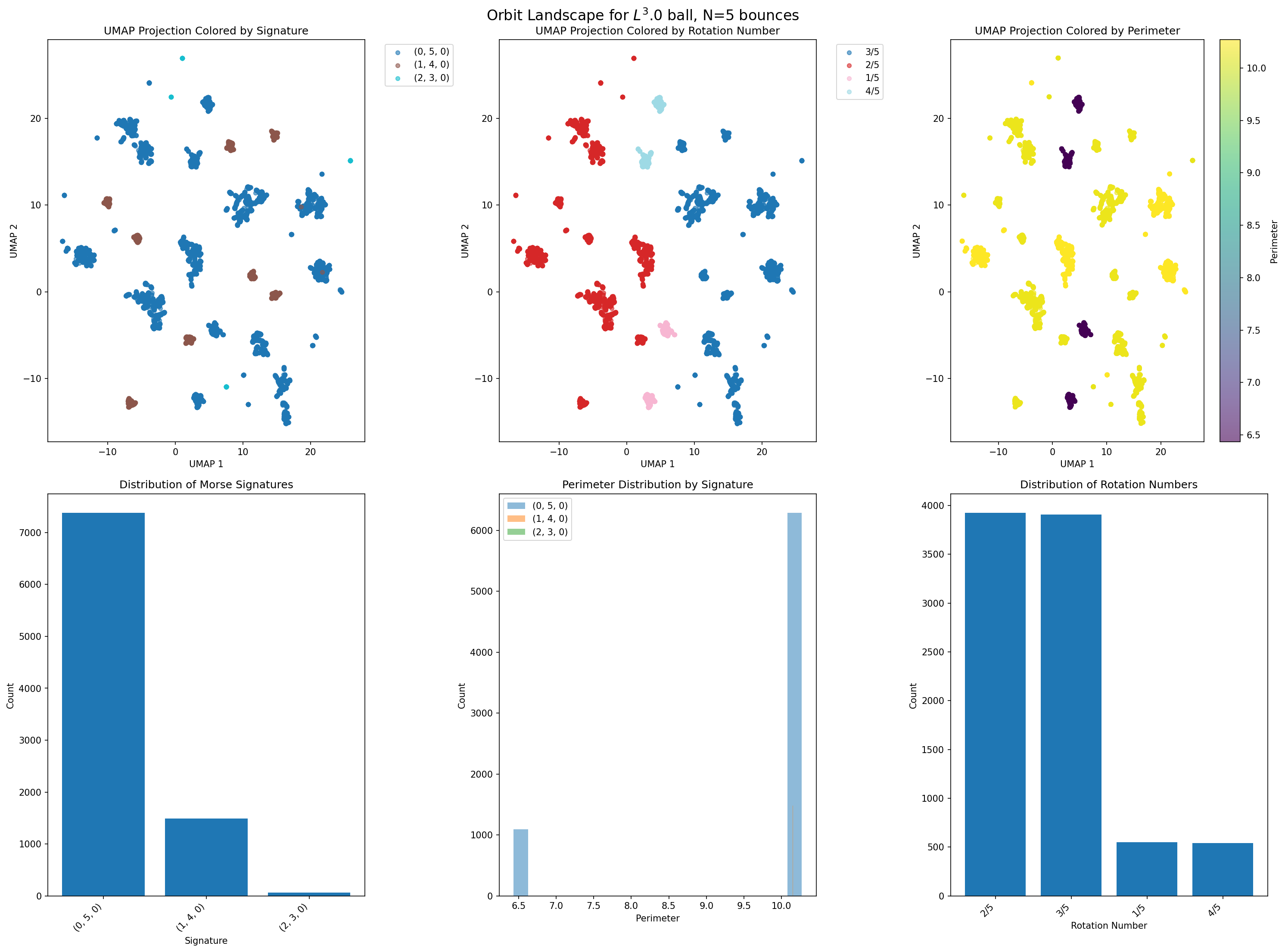}
\caption{UMAP visualization of 8,927 periodic orbits colored by: (a) Morse signature, (b) rotation number, and (c) perimeter. The projection reveals hierarchical clustering: orbits first separate by signature, then by rotation number within each signature, with continuous variation in perimeter creating gradients within clusters.}
\label{fig:umap}
\end{figure}

Key insights from the UMAP analysis (Figure~\ref{fig:umap}):

\begin{enumerate}
\item \textbf{Hierarchical structure}: The primary clusters correspond to Morse signatures, with sub-clusters for rotation numbers.
\item \textbf{Continuous families}: Within each signature-rotation combination, orbits form continuous manifolds in parameter space.
\item \textbf{Perimeter gradients}: The perimeter varies smoothly within clusters, suggesting continuous deformation families.
\item \textbf{Signature isolation}: Different signatures occupy distinct regions of the reduced space, indicating fundamentally different orbit geometries.
\end{enumerate}

\subsection{Comparison with Small-Scale Results}

Our initial experiments with 1,000 seeds found 149 orbits. The large-scale exploration reveals this was sampling only $\sim 1.7\%$ of the accessible orbit landscape at $N=5$. The power-law growth suggests even more orbits remain to be discovered.

\begin{figure}[h]
\centering
\includegraphics[width=0.8\textwidth]{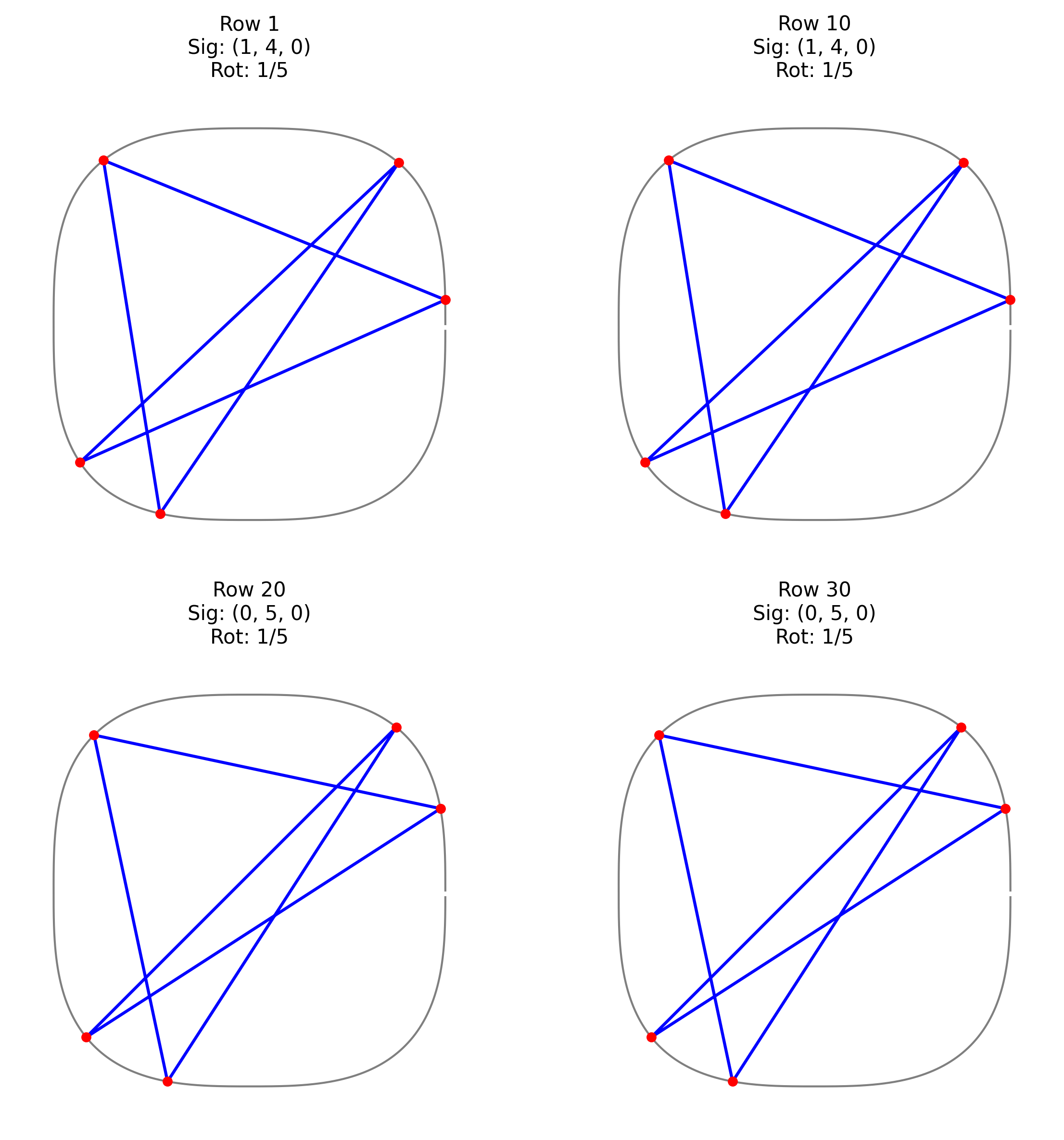}
\caption{A selection of periodic billiard orbits with \texorpdfstring{$N=5$}{N=5} bounces in the \texorpdfstring{$L^3$}{L3} ball, showing various rotation numbers (\texorpdfstring{$1/5$}{1/5}, \texorpdfstring{$2/5$}{2/5}, \texorpdfstring{$3/5$}{3/5}, \texorpdfstring{$4/5$}{4/5}) and Morse signatures. Each orbit is labeled with its signature \texorpdfstring{$(n_+, n_-, n_0)$}{(n+, n-, n0)} and rotation number \texorpdfstring{$r/s$}{r/s}. Star patterns (rotation numbers \texorpdfstring{$2/5$}{2/5}, \texorpdfstring{$3/5$}{3/5}) are prominently represented.}
\label{fig:grid_orbits}
\end{figure}

Figure~\ref{fig:grid_orbits} illustrates a small sample of the periodic orbits found by our method. While Birkhoff's theorem guarantees at least two orbits for each rotation number, our large-scale exploration found thousands of distinct orbits. Each corresponds to a different critical point of the perimeter functional, and they organize into a rich hierarchical structure revealed by dimensionality reduction techniques.

\section{Software Tools}

We provide three main utilities, available at \url{https://github.com/igorrivin/billiards}:

\begin{description}
\item[\texttt{orbit\_finder.py}] Main computation engine using JAX\cite{jax2018github} for finding orbits
\item[\texttt{verify\_orbit.py}] Standalone verification tool using pure NumPy
\item[\texttt{plot\_orbit.py}] Visualization tool with flexible input options
\end{description}

Example usage:
\begin{verbatim}
# Find orbits
python orbit_finder.py 3.0 5 1000

# Verify a specific orbit
python verify_orbit.py 3.0 "[0.0657, 0.375, 0.6843]"

# Visualize results
python plot_orbit.py --csv p3.0_N5_orbits.csv --row 1
\end{verbatim}

\section{Open Questions}

\begin{enumerate}
\item Can the observed signature patterns be explained by Morse theory or Lyusternik-Schnirelman theory?
\item Why do saddle orbits often have longer perimeters than local maxima?
\item What is the complete classification of possible signatures for given $(p, N)$?
\item How do the patterns change for $p \to 2^+$ (approaching the circle) or $p \to \infty$ (approaching the square)?
\item What explains the power-law growth $N_{\text{orbits}} \propto N_{\text{seeds}}^{0.714}$? Is there a theoretical bound on the number of periodic orbits?
\item Can the hierarchical clustering structure revealed by UMAP be explained by the underlying symplectic geometry?
\end{enumerate}

\section{Conclusion}

We have presented an efficient computational method for finding periodic billiard orbits in $L^p$ balls, with rigorous verification via Smale's criterion. Our large-scale exploration discovered nearly 9,000 distinct orbits for $N=5$ alone, far exceeding theoretical guarantees. The observed patterns in signatures and rotation numbers, combined with the hierarchical clustering revealed by UMAP, suggest a rich mathematical structure connecting dynamics, topology, and geometry. The power-law growth and lack of saturation indicate that the landscape of periodic billiard orbits is even richer than previously imagined, opening new avenues for both computational and theoretical investigation.

\appendix

\section{Algorithm Details}

\begin{algorithm}
\caption{Main Newton-Smale Algorithm}
\begin{algorithmic}
\REQUIRE Initial seeds $\{\theta^{(0)}_i\}_{i=1}^m$, shape parameter $p$, steps $k$
\ENSURE Certified orbits with signatures
\FOR{each seed $\theta^{(0)}$}
    \FOR{$j = 1$ to $k$}
        \STATE Compute $g = \nabla L(\theta)$, $H = \nabla^2 L(\theta)$
        \STATE Solve $Hd = -g$ for Newton direction $d$
        \STATE Update $\theta \leftarrow \theta + d$
    \ENDFOR
    \STATE Compute $\alpha = \beta \gamma$ using final $(g, H)$
    \IF{$\alpha < 0.15767$}
        \STATE Add to certified orbits with signature
    \ENDIF
\ENDFOR
\STATE Coalesce near-duplicates
\RETURN Unique certified orbits
\end{algorithmic}
\end{algorithm}
\bibliographystyle{alpha}
\bibliography{references}
\end{document}